# SOVIET PIONEERS OF FRACTIONAL CALCULUS AND ITS APPLICATION:
# III. TIMOTHY D. SHERMERGOR

Olga G. Novozhenova

Mech. Engin. Research Institute
Russian Academy of Sciences
M. Kharitonievsky Per. 4
101990- Moscow, RUSSIA

**Abstract:** In the first part of this survey we give a brief outline of T. D. Shermergor's biography and consider 3 of his works on the use of confluent hypergeometric function of the first kind (**CHGF)** for description asymmetric relaxation spectrum . The second part is devoted to Soviet works of the past (20th) that are in tune with this topic.
    This is the third part of the historical survey [8-9].
**AMS Subject Classification:**

## Part I.
### 1.Timothy D. Shermergor-Biography and 4 his works on CHGF

Timothy D. Shermergor was born in the city of Semipalatinsk, in a family of Dutch citizens who came to Russia in 1924 to help build socialism. His father, Dirk Shermergorn, was a construction specialist , mother, Franciska Shermergorn was a stomatologist and worked as an interpreter. Timothy later had a sister Anna (1930) and a brother Ilya (1937).

Dirk Shermergorn worked on the construction of Turksib, the railroad through Central Asia, and later on the construction of the subway in Moscow. From 1931 to 1938 the family lived in Moscow. In 1936, his father was arrested and, as an enemy of the people, was shot in november 1937. In 1938, Franciska Shermergorn, as the wife of the enemy of the people, was convicted and sent to a prison camp in the Mordovian region, along with her youngest son, Ilya. She died in 1953, Ilya Shermergorn, who became a professor at the Kazan Chemical Technology Institute, in 1993 was the first director of an open company NAPOR. Timothy and Anna (a retired doctor living in Chaikovsk, Perm region) were sent to various orphanages. At registration of Timothy in the orphanage, the last letter "n" of the surname was overlooked, and he was recorded as Timofey Dmitrievich Shermergor. Also, the date of birth was incorrectly mentioned, which later appeared in all his documents : 22.10.1929.

From 1938 to 1941 Timothy was in the orphanage was in the Kharkiv region, Ukraine. After the outbreak of the war with Germany, the orphanage was evacuated to town of Ocher, the Perm region, Urals. In 1946, after graduation at Shool, Timothy was awarded with a medal "For special achievements in education" and enrolled in the Faculty of Physics of Perm University.

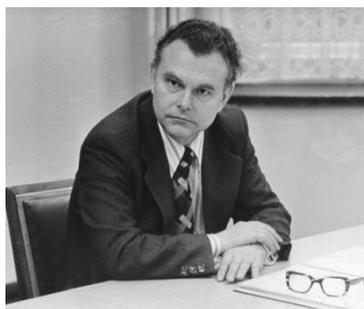

Fig.1: Timofey D.Shermergor
*(08.10.1928-18.07.1998)*

In 1951, he graduated with honors from the University on a speciality "Theoretical physics" and was sent to work at the Siberian Metallurgical Institute, Novokuznetsk, Kemerovo Region. Here, at the Department of Physics (1951-1960) Timofey Dmitrievich occupied posts of an assistant, a senior lecturer and assistant professor. In 1956 at Tomsk University (without training in postgraduate study) he defended his thesis for a candidate of science. In 1953 he married Uliana Aleksandrovna Timofeeva (1928 -1975), sons Yury (1955) and Eugene (1960).

In 1960, Timofey Dmitrievich moved to the city of Voronezh, worked as an associate professor of the Department of Physics of Voronezh

Institute of Technology. Late he headed the Department of Physics at the Voronezh Polytechnic Institute.

In 1966 Timofey Dmitrievich moved to work at MIET (Moscow Institute of Electronic Technology, Moscow, Zelenograd). Here he served as head of the Department of General Physics, then head of the Department of Theoretical Physics.

In 1966 he defended his doctoral thesis on the topic "Problems of the theory of mechanical relaxation in solid bodies ".

**1.** In 1967, T.D. Shermergor, in his work "The Rheological Characteristics of Viscoelastic Materials Having an Asymmetric Relaxation Spectrum" [10], for the first time justified the use of a confluent hypergeometric functions of the first kind to describe some properties of elastic- viscous bodies. It was shown that the fractional-exponential Rabotnov' kernels lead to a *symmetrical* bell-shaped distribution function of the logarithms of the relaxation times $H(\tau)$, although often the real functions $H(\tau)$ turn out to be essentially asymmetric[1] in the region of the transition of polymers from the glassy state to a high-elastic.
The asymmetric distribution function of logarithms of relaxation times corresponding to the Rzhanitsyn' kernel has a gentle decline in the direction of decreasing $\tau$ and a sharp decline in the direction of increasing. This contradicts experimental data [1]. It was proved that the correct asymmetry gives the use of **$_1F_1(α,1,x)$ that is a confluent hypergeometric function.** It can be defined by the series

$$_1F_1(\alpha,1,x) = 1 + \frac{\alpha}{(1!)^2}x + \frac{\alpha(\alpha+1)}{(2!)^2}x^2 + ... \frac{\alpha(\alpha+1)...(\alpha+n-1)}{(n!)^2}x^n + ...$$

For $\alpha = 1$, the function S (t) degenerates into an exponential function $_1F_1(1,1,-t/\tau) = \exp(-t/\tau)$.

This reproduces the transition to the standard linear body model. The image of the operator $P^*$ has the form

$$LP_\nu(p) = \frac{1+q_\nu}{(1+1/p\tau_\nu)^\alpha + q_\nu} \quad (\nu = 3,4),$$

where $q_3 = 1/m - 1$, $q_4 = m - 1$, $\tau_3 = \tau_\varepsilon$, $\tau_4 = \tau_\sigma$.

The fraction can be represented as a series, which in the case $|q_\nu| < 1|$, has the form 
$$LP_\nu(p) = (1+q_\nu)\sum_{n=0}^{\infty}(-q_\nu)^n(1+\frac{1}{p\tau_\nu})^{-\alpha(n+1)}$$

For $|q_\nu| > 1$ we obtain another series

$$LP_\nu(p) = (1 + \frac{1}{q_\nu})\sum_{n=0}^{\infty}(-q_\nu)^{-n}(1 + \frac{1}{p\tau_\nu})^{\alpha n}$$

The inversion of the expressions for $LP_\nu(p)$ gives the following

$$P_\nu(t) = (1 + q_\nu)\sum_{0}^{\infty}(-q_\nu)^n \,_1F_1[\alpha(n+1), 1, -t/\tau_\nu]$$

$$P_\nu(t) = (1 + \frac{1}{q_\nu})\sum_{0}^{\infty}(-q_\nu)^{-n} \,_1F_1[-\alpha n, 1, -t/\tau_\nu]$$

At the end of the article, the gratitude is expressed to Yu. N. Rabotnov for the discussion.

**2.** The results of further research on the topic [10] were published in co-authorship with Vladimir Listovnichy [4-6].

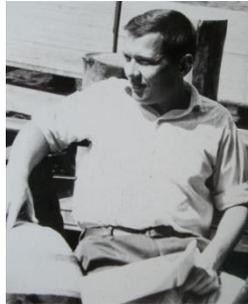

*Fig. 2: V. Listnovnichy (born 06.07.1940)*

*Vladimir F. Listnovnichy studied at the postgraduate school of MIET (Moscow Institute of Electronic Technology, 1966-1969) in 1971, defended his thesis on "Some problems of mechanical and dielectric relaxation in inhomogeneous materials"[4] and received a degree of candidate of physical and mathematical science. After graduate School and obtaining a scientific degree, V.F.Listovnichy left MIET.*

The creep of viscoelastic bodies were considered in [5] with kernel in the form of a confluent hypergeometric function. The integral operators of the elastic moduli and the compliances $J^*$ are defined by the equalities

$$M^* = V_\infty(1 - qR^*) = M_0(1 + pS^*D), \quad D = \partial/\partial t$$

A kernel in the form of a confluent hypergeometric function can be defined by the following expression for the Fourier transform of the integral operator of the elastic modulus

$$M(i\omega) = M_0 + (M_\infty - M_0)(1 + 1/i\omega\tau_\varepsilon)^{-\alpha}$$

where $\tau_\varepsilon$ is the characteristic relaxation time, α is the fractional parameter. It follows that

$$S(t) = {}_1F_1(\alpha, 1, -\theta), \quad R(t) = \frac{\alpha}{\tau_\varepsilon} {}_1F_1(1+\alpha, 2, -\theta), \quad \theta = \frac{t}{\tau_\varepsilon}.$$

An analogous approach was also used *for the fractional-exponential kernels of Yu. N. Rabotnov,* with folowing results

$$M(i\omega) = M_\infty - \frac{\Delta M}{1 + (i\omega\tau_\varepsilon)^\alpha}, \qquad J(i\omega) = J_\infty - \frac{\Delta J}{1 + (i\omega\tau_\sigma)^\alpha}$$

$$H(\tau) = \frac{1}{2\pi} \frac{\sin\alpha\pi}{ch[\alpha \ln(\tau/\tau_\varepsilon)] + \cos\alpha\pi}, \quad L(\tau) = \frac{1}{2\pi} \frac{\sin\alpha\pi}{ch[\alpha \ln(\tau/\tau_\sigma)] + \cos\alpha\pi}$$

$$\varphi(t) = I_\alpha(\theta_\varepsilon), \quad \psi(t) = 1 - I_\alpha(\theta_\sigma)$$

$$I_\alpha(\theta) = 1 - (eh)_\alpha * 1 = \frac{\sin\alpha\pi}{\pi} \int_0^\infty x^{\alpha-1}(1 + x^{2\alpha} + 2x^\alpha \cos\alpha\pi)^{-1} e^{-\theta x} dx \quad (A)$$

$$\theta_\varepsilon = \frac{t}{\tau_\varepsilon}, \quad \theta_\sigma = \frac{t}{\tau_\sigma}, \left(\frac{\tau_\varepsilon}{\tau_\sigma}\right)^\alpha = \frac{M_0}{M_\infty}, \quad 0 < \alpha < 1$$

In contrast to the previous case, the distribution functions $H(\tau)$ and $L(\tau)$ in logarithmic coordinates were symmetric and only shifted one relative to the other.

**3.** *Relaxation of stress and creep of some linear viscoelastic solids [6]*

In this paper a model of Rabotnov's (*Shermergor 's in the manuscript*) viscoelastic solid was constructed, which takes into account the asymmetry of the relaxation spectrum and have been introduced the **splitting property** of the kernel. Resolvent kernels are defined in an explicit form both in terms of a series and in the integral representation. A detailed comparison is performed for the fractional-exponential kernels, for the Rzhanitsyn-Davidson' kernels and for the kernel in the form of a confluent hypergeometric function and the asymptotic ratios for large time values are obtained.

The simplest kernels that take into account the asymmetry of the relaxation spectrum are the Rzhanitsyn-Davidson' kernels and the confluent hypergeometric function $_1F_1(\alpha+1,2,\theta)$, where $0<\alpha<1$. The disadvantage of these kernels is that they are not possess the splitting property. In connection with this, it is of interest to construct models of such viscoelastic solid, that take into account the asymmetry of the spectrum and lead to splittable operators of both the elastic moduli $M^*$ and the compliancy $J^*$, analogously to the one proposed by Yu. N. Rabotnov for the Abel's kernel.

They proceed from the following concepts for complex compliances of simple hereditary solid:

$$J(i\omega) = J_\infty + J_\infty(i\omega\tau)^{-\alpha} = J_\infty[1+\kappa'(i\omega)^{-\alpha}] \tag{1}$$

$$J(i\omega) = J_\infty + \Delta J(1+i\omega\tau)^{-\alpha} = J_\infty[1+\kappa''(1/\tau+i\omega)^{-\alpha}] \tag{2}$$

$$J(i\omega) = J_0 - \Delta J(1+1/i\omega\tau)^{-\alpha} = J_\infty\{1+\kappa'''[\tau^{-\alpha}-(\tau+1/i\omega)^{-\alpha}]^{-1}\} \tag{3}$$

$$\kappa' = \tau^{-\alpha}, \quad \kappa'' = (1/m-1)\tau^{-\alpha}, \quad \kappa''' = (1-m-1)\tau^{\alpha}, \quad m = J_\infty/J_0 \tag{4}$$

Here $J_0$ and $J_\infty$ are, respectively, relaxed and non-relaxed compliance, $\Delta J$ is their difference, ω is a cyclic frequency, τ is the characteristic relaxation time, and α-fractional index.

The transition in the expressions (1) – (3) from the image leads to the Abel's, Rzhanitsyn-Davidson's kernels and the kernels in the form of the confluent hypergeometric function type.

The resolvent operator $R^*$ is defined by

$$R^*(\lambda,T) = (\lambda E^* + T^*)^{-1}, \tag{5}$$

where $E^*$ is the unit operator, $T^*$ is the basic operator and $\lambda$ is a number. From (5) it follows that the Hilbert identity (the splitting condition)

$$R^*(\lambda,T)R^*(\mu,T) = \frac{1}{\mu-\lambda}\left[R^*(\lambda,T) - R^*(\mu,T)\right] \tag{6}$$

and the rule of lowering the degree of the operator is given below

$$R^{*2}(\lambda,T) = -\frac{\partial}{\partial\lambda}R^*(\lambda,T) \tag{7}$$

The operators of shear modulus of elasticity and compliance are expressed in terms of the resolvent operator by means of relations

$$M^* = M_\infty[1-\kappa R^*(\lambda,T)], \quad J^* = J_\infty[1+\kappa R^*(\lambda-\kappa,T)], \tag{8}$$

$$J^* = J_\infty[1+\kappa R^*(\mu,T)], \quad M^* = M_\infty[1-\kappa R^*(\mu+\kappa,T)], \tag{9}$$

where the second pair is obtained from the first by replacing $\lambda - \kappa = \mu$. In both equations, the first place is given to direct operators that determine the hereditary properties of the medium through the resolvent operator $R^*$, and the second place is occupied by the inverse operators, which are expressed through the same resolvent operator by means of a parameter shift.

They require that for $\lambda = 0$, the resolvent operator (5) leads to one of the equalities (1) - (3). This gives the following result

$$LT_1 = (i\omega)^\alpha \equiv LI_{-\alpha} \tag{10}$$

$$LT_2 = (1/\tau + i\omega)^\alpha \equiv LK_{-\alpha} \tag{11}$$

$$LT_3 = [\tau^{-\alpha} - (\tau + 1/i\omega)^{-\alpha}]^{-1} \equiv (LN_\alpha)_{-1} \tag{12}$$

Here L is the operator of the integral Fourier transform.
Equations (5), (8),(9) and one of the formulas (10) - (12) completely determine **the hereditary properties of the Rabotnov's (*Shermergor's*) viscoelastic solid.**
Using designations
$$R(\lambda, I_{-\alpha}) = (eh)_\alpha(\lambda), \quad R(\lambda, K_{-\alpha}) = Q_\alpha(\lambda), \quad R(\lambda, N_\alpha) = P_\alpha(\lambda),$$
operators of elastic moduli and compliances can be written in the standard form

$$M^* = M_\infty[1 - \kappa(eh)_\alpha^*(\lambda)], \quad J^* = J_\infty[1 + \kappa(eh)_\alpha^*(\mu)] \tag{13}$$

$$M^* = M_\infty[1 - \kappa Q_\alpha^*(\lambda)], \quad J^* = J_\infty[1 + \kappa Q_\alpha^*(\mu)] \tag{14}$$

$$M^* = M_\infty[1 - \kappa P_\alpha^*(\lambda)], \quad J^* = J_\infty[1 + \kappa P_\alpha^*(\mu)], \tag{15}$$

with

$$(eh)_\alpha(\lambda, t) = t^{\alpha-1} \sum_0^\infty \frac{(-1)^n (t/\tau)^{\alpha n}}{\Gamma[\alpha(n+1)]},$$

$$Q_\alpha(\lambda, t) = \exp(-\frac{t}{\tau})(eh)_\alpha(\lambda, t),$$

$$P_\alpha(\lambda, t) = \frac{1}{\tau^{\alpha+1}} \sum_{k=0}^\infty \frac{n_\varepsilon^k}{(n_\varepsilon + 1)^{k+1}} \{\alpha(k+1)\,_1F_1\left[\alpha(k+1), 2, -t/\tau\right] -$$
$$-k\alpha\,_1F_1(k\alpha + 1, 2, -t/\tau)\},$$

where $_1F_1(\alpha, c, x)$ - is a confluent hypergeometric function, and for basic operators (10)-(12) the parameter $\lambda$ is defined accordingly to the relations

$$\lambda = \mu + \kappa = \tau_\varepsilon^{-\alpha}, \quad \kappa = \lambda(1-m) = \mu(1/m-1)$$

$$\lambda - \kappa = \mu = \tau_\sigma^{-\alpha}, \quad m = \tau_\varepsilon^\alpha / \tau_\sigma^\alpha$$

$$\lambda = n_\varepsilon \tau^{-\alpha}, \quad \mu = n_\sigma \tau^{-\alpha}, \quad m = (1+n_\sigma)/(1+n_\varepsilon)$$

$$\kappa = (1-m)(n_\varepsilon + 1)\tau^{-\alpha} = (1/m-1)(n_\sigma + 1)\tau^{-\alpha}$$

$$\lambda = n_\varepsilon \tau^\alpha, \quad \mu = n_\sigma \tau^\alpha, \quad m = (1+n_\sigma)/(1+n_\varepsilon)$$

$$\kappa = (1-m)(n_\varepsilon + 1)\tau^\alpha = (1/m-1)(n_\sigma + 1)\tau^\alpha$$

Earlier [4], a method was proposed for finding integral representations of resolvent kernels and the formula was obtained for convolution of a fractional exponential function with unity. They use the same method to calculate the convolutions of operators $Q_\alpha^*$ and $P_\alpha^*$ with the unit. The corresponding calculation gives the following results

$$\tau^{-\alpha} Q_\alpha^*(\lambda) \cdot 1 = \frac{1}{n_\varepsilon + 1} - \frac{\sin \alpha\pi}{\pi} \int_1^\infty \frac{\xi}{\xi^2 + 2\xi n_\varepsilon \cos\alpha\pi + n_\varepsilon^2} \exp(-\theta x) \frac{dx}{x} \quad (16)$$

$$\xi \equiv (x-1)^\alpha, \quad n_\varepsilon = \lambda \tau^\alpha, \quad \theta = t/\tau$$

$$\tau^\alpha P_\alpha^*(\lambda) \cdot 1 = \frac{1}{n_\varepsilon + 1} - \frac{\sin \alpha\pi}{\pi} \int_0^1 \frac{\xi}{\xi^2(1+n_\varepsilon)^2 - 2\xi(1+n_\varepsilon)n_\varepsilon \cos\alpha\pi + n_\varepsilon^2} \exp(-\theta x) \frac{dx}{x} \quad (17)$$

$$\xi \equiv (1/x-1)^\alpha, \quad n_\varepsilon = \lambda \tau^{-\alpha}$$

As an example, it was considered the relaxation of stresses in an unbounded solid in which there is a cavity of radius R.
For all three cases of solid models (13) - (15), integral representations allow us to obtain *asymptotic estimates for large times*. For sufficiently large t and small τ, it follows from equalities (16) and (17) that

$$\kappa Q_\alpha^*(\lambda - \kappa_1) \cdot 1 \approx k(1-m)\left[1 - \frac{k\alpha(1+\lambda_0 - k)^{-2}}{\Gamma(1-\alpha)} \theta^{-\alpha-1} e^{-\theta}\right]$$

$$\kappa P_\alpha^*(\lambda - \kappa_1) \cdot 1 \approx k(1-m)\left[1 - \frac{k\theta^{-\alpha}}{\Gamma(\alpha-1)(\lambda_0 + 1)}\right]$$

$$\kappa(eh)^*_\alpha(\lambda-\kappa_1)\cdot 1 \approx k(1-m)\left(1-\frac{k\theta^{-\alpha}}{\Gamma(1-\alpha)}\right)$$

$$k^{-1} \equiv 1-4/5\alpha_1(1-m), \quad \theta = t/\tau$$

At the end of the article, the gratitude is expressed to Yu. N. Rabotnov for the discussion.

***N.B.*** *All articles[5-6,10], with the permission of Yu.T. Shermergor and V.F. Listovnichy, available on the site of IMASH (Mechanical Engineering Research Institute of Russian Academy of Sciences).*

**4.** T.D. Shermergor was published in 1977 a book [11], in the Preface again gratitude was expressed to Yu.Rabotnov for the discussion of the manuscript. In §1 of Chapter 8 he cited the results of [10] and [5-6] (references 70, 71 and 153), dealing only with the algebra of Rabotnov's eh-functions , without mentioning exotic confluent hypergeometric functions of the first kind (Kummer functions). Nowadays, half a century later, the Kummer functions are included in the standard set of MatLab, Maple, MathCad, Wolfram Mathematica functions.

Later, T. D. Shermergor, head of the Department of physics of MIET, was mainly engaged in the tasks of piezoelectricity and the publication of monographs of courses for students [2-3, 12-14].

## 2. The further research on asymmetric relaxation spectrum and confluent hypergeometric functions

The results, published in articles [5-6] and [10], found little application in the works of contemporaries. E.L. Sinaisky, M.I. Rozovsky's follower [8,17], considering the problem of realizing the function of a hereditary operator, acting upon a certain function of time, used the Laplace transform for operators with Yu.N.Rabotnov 's and A.R. Rzhanitsyn's kernels . He obtained the formulas, that reduce the given problem, to the calculation of the quadrature. For the fractional-exponential Rabotnov's kernels, the solution obtained in [5] was used

$$I_\alpha(\theta) = 1 - (eh)_\alpha * 1 = \frac{\sin \alpha\pi}{\pi} \int_0^\infty x^{\alpha-1}(1 + x^{2\alpha} + 2x^\alpha \cos \alpha\pi)^{-1} e^{-\theta x} dx.$$

M. A. Koltunov [9] was devoted a whole page of paper [10], noticed:
"To account for the asymmetry of the relaxation spectrum in the transition region polymers from a glass-like state to a highly elastic state, T.D. Shermergor suggests using confluent hyperbolic functions of the form

$$_1F_1(\alpha, 1, x) = 1 + \frac{\alpha}{(1!)^2} x + \frac{\alpha(\alpha+1)}{(2!)^2} x^2 + ... \frac{\alpha(\alpha+1)...(\alpha+n-1)}{(n!)^2} x^n + ...$$

The method of determining the kernel parameters from experimental data was not discussed».
*Note that in [10] the definition of α for polyisobutylene was given in two ways,* $\alpha = 0,61; \ \alpha = 0,63$.

*V.D. Strauss suggested* time integral characteristics of four parameter descriptions of relaxation processes with asymmetric spectrum[15-16]. He received the relationships between the kernel of relaxations and its resolvent, functions of the relaxation and creep function.

In the study of dielectrics, a four-parameter dispersion formula Havriliax-Negami is used, which allows to regulate the asymmetry of spectra in a wide range-from symmetric to explicitly asymmetric. The dispersion formula is a generalization of variances according to Debye, Coule-Coule and Davidson-Coule, special cases of it are also the models of mechanical relaxation Rabotnov and Rzhanitsyn.

For the frequency dependence of the complex module, the four-parameter Gavrilyak-Negami dispersion formula can be written as

$$\tilde{M}(j\omega) = M_\infty - \frac{\Delta M}{[1 + (j\omega\tau_0)^\alpha]^\beta}; \quad 0 < \alpha \leq 1; \quad 0 < \beta \leq 1; \quad (18)$$

$$\Delta M = M_\infty - M_0,$$

where $\alpha$, $\beta$ -parameters determining the shape of the spectrum, $\tau_0$ -the characteristic attenuation time of the process (the most likely time of macroscopic relaxation); $M_\infty$, $M_0$ –instantaneous and equilibrium modules; $j = \sqrt{-1}$. A function of relaxation rates $R(t)$ and normalized integrated module

$$\tilde{m}(j\omega) = \frac{\tilde{M}(j\omega) - M_\infty}{\Delta M} = \frac{1}{[1 + (j\omega\tau_0)^\alpha]^\beta} \tag{19}$$

are related by the inverse Laplace transform

$$R(t) = \frac{1}{2\pi j} \int_{c-j\infty}^{c+j\infty} \tilde{m}(p) \exp(pt) dp,$$

where the parameter $p = c + j\omega$.

By decomposing (19) in the power series and applying the transformation $\left(\frac{1}{p}\right)^{\nu+1} = \frac{t^\nu}{\Gamma(\nu+1)}$, to the relaxation kernel, corresponding to (18), was obtained

$$R(t) = \frac{1}{\tau_0 \Gamma(\beta)} \sum_{i=0}^{\infty} \frac{(-1)^i \Gamma(\beta+i) \left(\frac{t}{\tau_0}\right)^{\alpha(\beta+i)-1}}{\Gamma(i+1)\Gamma[\alpha(\beta+i)]}.$$

The corresponding relaxation function has the form

$$\int_0^t R(t)dt = \frac{1}{\Gamma(\beta)} \sum_{i=0}^{\infty} \frac{(-1)^i \Gamma(\beta+i) \left(\frac{t}{\tau_0}\right)^{\alpha(\beta+i)}}{\Gamma(i+1)\Gamma[\alpha(\beta+i)+1]}$$

To find the analytical form of the creep kernel $K(t)$, which is the resolvent of the function $R(t)$, the resolvent ratio of Volterra is used

$$\tilde{K}(p) - \tilde{R}(p) = \tilde{K}(p)\tilde{R}(p),$$

where $\tilde{K}(p)$ and $\tilde{R}(p)$ are the Laplace image of functions $K(t)$ and $R(t)$ respectively

The special cases are considered and the results are summarized in table

Kernels of $R(t)$ and $K(t)$ with four parameters

| Parameter value | | Kernel $R(t)$ | Resolvent $K(t)$ |
|---|---|---|---|
| $\alpha$ | $\beta$ | | |
| $0 < \alpha \leq 1$ | $0 < \beta \leq 1$ | $\dfrac{1}{\tau_0 \Gamma(\beta)} \sum\limits_{i=0}^{\infty} \times \dfrac{(-1)^i \Gamma(\beta+i)(t/\tau_0)^{\alpha(\beta+i)-1}}{\Gamma(i+1)\Gamma[\alpha(\beta+i)]}$ | $\dfrac{1}{\tau_0} \sum\limits_{i=0}^{\infty} \sum\limits_{n=1}^{\infty} \times \dfrac{(-1)^i \Gamma(n\beta+i)(t/\tau_0)^{\alpha(n\beta+i)-1}}{\Gamma(n\beta)\Gamma(i+1)\Gamma[\alpha(n\beta+i)]}$ |
| $0 < \alpha \leq 1$ | 1 | $\dfrac{1}{\tau_0} \sum\limits_{i=0}^{\infty} \dfrac{(-1)^i (t/\tau_0)^{\alpha(1+i)-1}}{\Gamma[\alpha(1+i)]}$ <br> (Rabotnov, Cole-Cole) | $\dfrac{(t/\tau_0)^{\alpha-1}}{\tau_0 \Gamma(\alpha)}$ <br> [Abel(Duffing), Curie] |
| 1 | $0 < \beta \leq 1$ | $\dfrac{1}{\tau_0 \Gamma(\beta)} \left(\dfrac{t}{\tau_0}\right)^{-(1-\beta)} \times \exp\left(-\dfrac{t}{\tau_0}\right)$ <br> (Rzhanitsyn, Davidson-Cole) | $\dfrac{\exp(t/\tau_0)}{t} \sum\limits_{n=1}^{\infty} \dfrac{(t/\tau_0)^{n\beta}}{\Gamma(n\beta)}$ <br> (Koltunov) |
| 1 | 1 | $\dfrac{1}{\tau_0} \times \exp\left(-\dfrac{t}{\tau_0}\right)$ <br> (elementary relaxation process, Debye) | ——— |

In [16] a method for calculating the kernels of heredity and integrals from them was given, based on the application of the inverse Laplace transform to the corresponding functions of their images.

G.A. Vanin, after the death of Yu.N.Rabotnov, became the head of the laboratory of Composition Materials in IMASh RAS. His attention was drawn to the possibility of an analytical description of the viscoelastic properties of a structured material in terms of their distribution functions in combination with the use of Kummer functions satisfying the determining equations [20-22].

To describe asymmetrically distributed positive random parameters (diameters of circular or spherical inclusions in heterogeneous media, distances between their centers, etc.), new distribution functions are proposed [20], their central moments are determined. The distribution density was entered

$$p(x) = Ax^b \exp(-\frac{x^2}{2\sigma^2}) sh(\frac{\alpha x}{\sigma}), \quad (0 \le x < \infty).$$

Here a, b-dimensionless parameters, established, as well as $\sigma$, on the basis of experimental data processing, A-normalizing factor, determined from the condition

$$\int_0^\infty p(x)dx = 1, \quad A = [2^{b/2}\alpha\sigma^{1+b}\Gamma(1+\tfrac{b}{2})_1F_1(1+\tfrac{b}{2},\tfrac{3}{2},\tfrac{a^2}{2})]^{-1},$$

where $_1F_1(1+\tfrac{b}{2},\tfrac{3}{2},\tfrac{a^2}{2})$ -is a confluent hypergeometric function.

In [21-22], histograms and density of distributions were constructed from data on the distribution of microspheres diameters, distances between them, angles between the lines drawn through the centers of neighboring microspheres. The formula for the initial moments was obtained.

A quarter of a century later [10], in the work of Yu. V. Suvorova, [18-19], Yu.N.Rabotnov's follower, a *hypergeometric function* appears with regard to the inversion of the relationship between the intensities of stresses and strains in a nonlinear form

$$\varphi(\varepsilon_i) = (1+K^*)\sigma_i.$$

The obtained dependence is similar to equation Lederman-Rozovsky, but with the same functions of nonlinearity; for $K(t) = kt^{-\alpha}$ we have relationship expressed in terms of gamma function and Rabotnov' eh-function:

$$\sigma_i = \varphi(\varepsilon_i) - \int_0^t k\Gamma(1-\alpha)(\text{eh})_\alpha(-k, t-\tau)\varphi[\varepsilon_i(\tau)]d\tau$$

For $\varphi(\varepsilon_i) = a\ln(1+b\varepsilon_i)$, $\varepsilon_i = \dot{\varepsilon}_i t$ and

$$(\text{eh})_\alpha(-k,t) = t^{-\alpha}\sum_{n=0}^\infty \frac{[-k\Gamma(1-\alpha)]^n t^{n(1-\alpha)}}{\Gamma[(n+1)(1-\alpha)]}$$

received

$$\sigma_i = \varphi(\varepsilon_i) + \sum_{n=1}^{\infty} \frac{a[-k\Gamma(1-\alpha)]^n b\dot{\varepsilon}_i t^{n(1-n)+1}}{\Gamma[(n+1)(1-\alpha)]} \times {}_2F_1(1,1;2+n(1-\alpha);-b\varepsilon_i).$$

When $n \to \infty$ the series of the hypergeometric function ${}_2F_1(1,1;2+n(1-\alpha);-b\varepsilon_i)$ tends to 1, and the whole series tends to 0. From the condition $\alpha \approx 1$ we have $n(1-\alpha)+2 \approx 2$ and the dependence is obtained ${}_2F_1(1,1;2;-b\varepsilon_i) = \frac{\ln(1+b\varepsilon_i)}{b\varepsilon_i}$, used to process the experiment [16-17].

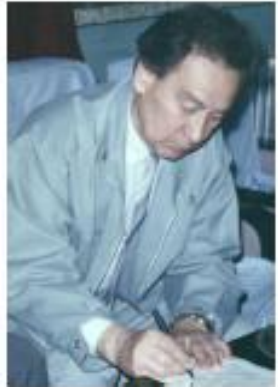

Fig.3. Jeorge A.Vanin
(1930-2008)

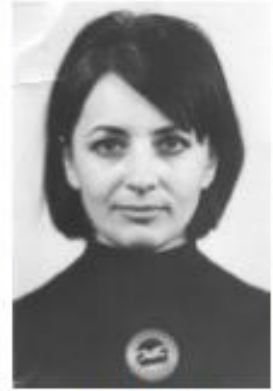

Fig.4 Julia V.Suvorova
(1939-2009)

### Acknowledgements

*The author expresses his deep gratitude to Yuri Timofeevich Shermergor, the son of T.D.Shermergor, for the photographs and brief biography of his father. Special thanks Vladimir F. Listovnichy, for providing a photo and the brief biography.*